\documentclass[12pt]{amsart}
\usepackage{amsmath, amsthm, amscd, amsfonts}
\newtheorem{Lemma}{Lemma}

 \newfont{\Bbbb}{msbm10 scaled\magstephalf}
 \def\be{\begin{equation}}
 \def\ee{\end{equation}}

\def\bn{\begin{equation}}
\def\en{\end{equation}}
\def\br{\begin{center}}
\def\er{\end{center}}
\def\by{\begin{array}}
\def\ey{\end{array}}

\def\begy{\begin{eqnarray}}
\def\endy{\end{eqnarray}}
\def\bey*{\begin{eqnarray*}}
\def\eny*{\end{eqnarray*}}
\def\ber{\begin{tabular}}
\def\enr{\end{tabular}}

\def\bt{\begin{flushright}}
\def\et{\end{flushright}}

 \def\bea{\begin{equation}\begin{array}{ll}}
 \def\eea{\end{array}\end{equation}}

\begin{document}
\title[Difference of composition operators]{Difference of composition operators  in the Polydiscs}
\author[Z.S. Fang and Z.H.Zhou]{Zhongshan Fang \and Zehua Zhou$^*$ }
\address{\newline Department of Mathematics\newline
Tianjin Polytechnic University
\newline Tianjin 300160\newline P.R. China.}

\email{fangzhongshan@yahoo.com.cn}

\address{\newline Department of Mathematics\newline
Tianjin University
\newline Tianjin 300072\newline P.R. China.}
\email{zehuazhou2003@yahoo.com.cn}



\date{}
\thanks{\noindent $^*$Zehua Zhou, corresponding author. Supported in part by the National Natural Science Foundation of
China (Grand Nos. 10671141, 10371091).}

\begin{abstract}
This paper gives some simple estimates of the essential norm for the
difference of composition operators induced by $\varphi$ and $\psi$
acting on bounded function space in the unit polydiscs $U^n$, where
$\varphi(z)$ and $\psi(z)$be holomorphic self-maps of $U^n$. As its
applications, a characterization of compact difference is given for
composition operators acting on the bounded function spaces.
\end{abstract}

\maketitle


\section{Introduction}

Let $U^n$ be the unit polydiscs of $C^n$ with boundary $\partial
U^n$.  The class of all holomorphic functions on domain $U^n$ will
be denoted by $H(U^n)$. Let
$\varphi(z)=(\varphi_1(z),\cdots,\varphi_n(z))$ and
$\psi(z)=(\psi_1(z),\cdots,\psi_n(z))$ be  holomorphic self-maps of
$U^n$. Composition operator is defined by:
$$
C_\varphi(f)(z)=f(\varphi(z))$$ for any $ f\in H(U^n)$ and $z\in
U^n$.

We recall that the essential norm of a continuous linear operator
$T$ is the distance from $T$ to the compact operators, that is,
$$\|T\|_e=\inf\{\|T-K\|: K \mbox{ is compact}\}.$$ Notice that
$\|T\|_e=0$ if and only if $T$ is compact, so that estimates on
$\|T\|_e$ lead to conditions for $T$ to be compact.

In the past few years, boundedness and compactness of composition
operators between several spaces of holomorphic functions have been
studied by many authors. Recently, there have been many papers
focused on studying the mapping properties of the difference of two
composition operators, i.e., an operator of the form
$$T=C_{\varphi}-C_{\psi}.$$

In \cite{Mac2}, MacCluer et al., characterize the compactness of the
difference of composition operator on $H^{\infty}$ spaces by
Poincar\'{e} distance. In \cite{Toe} and \cite{GorM}, Carl and
Gorkin et al., independently extended the results to
$H^{\infty}(B_n)$ spaces, they described compact difference by
Carath\'{e}odory psedo-distance on the ball, which is the
generalization of Poinar\'{e} distance. In \cite{Moo}, Moorhouse
shew that if the pseudo-hyperbloic distance between the image values
$\varphi$ and $\psi$ converges to zero as $z\to\zeta$ for every
point $\zeta$ at which $\varphi$ and $\psi$ have finite angular
derivative then the difference $C_{\varphi}-C_{\psi}$ yields a
compact operator. Building on this foundation, this paper gives some
simple estimates of the essential norm for the difference of
composition operators induced by $\varphi$ and $\psi$ acting on
bounded function space in the unit polydiscs $U^n$, where
$\varphi(z)$ and $\psi(z)$be holomorphic self-maps of $U^n$. As its
applications, a characterization of compact difference is given for
composition operators acting on the bounded function spaces by
Carath\'{e}odory distance on $U^n$.

\section{Notation and background}

Let $D$ be the unit disc in $C$, then the pseudo-hyperbolic distance
on $D$ is defined by $\beta(z,w)=|\frac{z-w}{1-\overline{z}w}|$ for
$z,w \in D$. By $U^n$ denote the unit polydiscs of $C^n$, and by
$H^{\infty}$ denote the Banach space of bounded analytic functions
on $U^n$ in the sup norm.

{\bf Definition 1} The Poincar\'{e} distance $\rho$ on $D$ is
$$ \rho(z,w):=\tanh^{-1}\beta(z,w)
$$ for $z,w \in D$.

{\bf Definition 2} The Carath\'{e}odory pseudo-distance on a domain
$G\subset C^n $ is given by
$$c_{G}(z,w):=\sup\{\rho(f(z),f(w)):f\in H(G,D)\}$$ for $z, w \in
G.$

If we put $c_{G}^{\ast}(z,w):= \sup\{\beta(f(z),f(w)):f\in H(G,D)\}$
for $z, w \in G$, it is clear that
$$c_G=\tanh^{-1}(c_{G}^{\ast})\geq c_{G}^{\ast}.$$
Next we introduce the following pseudo-distance on
G$$d_G(z,w):=\sup\{|f(z)-f(w)|:f\in H(G,D)\}.$$ For the case $G=D$,
it is easy to show that
$$d_D(z,w)=\frac{2-2\sqrt{1-\beta(z,w)^2}}{\beta(z,w)}.$$ So the
Poincar\'{e} metirc on $D$
$$\rho(z,w):=\tanh^{-1}\beta(z,w)=\log\frac{2+d_D(z,w)}{2-d_D(z,w)}.$$

Clearly for $z,w \in G,$
\begin{eqnarray*}
d_G(z,w) &=& \sup\{|g(f(z))-g(f(w))|:g\in H(D,D),f\in
H(G,D)\}\\
&=& \sup\limits_{f\in H(G,D)}d_{D}(f(z),f(w))
\end{eqnarray*}

Since the map $t\rightarrow \log\frac{2+t}{2-t}$ is strictly
increasing on $[0,2),$ it follows that
\begin{eqnarray*}
\log\frac{2+d_G}{2-d_G}&=&\sup\limits_{f\in
H(G,D)}\log\frac{2+d_D(f(z),f(w))}{2-d_D(f(z),f(w))}\\
&=& \sup\limits_{f\in H(G,D)}\rho(f(z),f(w))\\
&=& c_G(z,w)
\end{eqnarray*}
 or equivalently for any
domain $G$ and any $ z, w \in G.$
\begin{eqnarray*}
d_G(z,w)&=&\frac{2-2\sqrt{1-(\tanh c_G(z,w))^2}}{\tanh c_G(z,w)}\\
&=& \frac{2-2\sqrt{1-( c_{G}^{\ast}(z,w))^2}}{c_{G}^{\ast}(z,w)}
\end{eqnarray*}

It is well known that $c_{U^n}^{\ast}(z,w)=\max\limits_{1\leq j\leq
n}{\beta(z_j,w_j)}.$ So we have
$$d_{U^n}(z,w)=\frac{2-2\sqrt{1-(\max\limits_{1\leq j\leq
n}\beta(z_j,w_j))^2}}{\max\limits_{1\leq j\leq n}\beta(z_j,w_j)}$$

Before proving the main theorem, we give first some symbol. For any
$0<\delta<1$, define $$E^j_{\delta}:=\{z\in U^n| \;\;
|\varphi_j(z)|\vee|\psi_j(z)|>1-\delta\},$$ and we put
$E_\delta=\cup_{j=1}^n \; E^j_{\delta},$where $\vee$ means the
maximum of two real numbers.

\begin{Lemma}  Let $\{z_n\}$ be a sequence in D with $|z_n|\rightarrow 1$ as $n\rightarrow \infty$, then there is a subsequence $\{z_{n_i}\} $,a number $M \geq 1$ and
a sequence of functions $f_m\in H^\infty(D) $ such that
\begin{eqnarray*}
&i)& \;\;\; f_m(z_{n_k})=\delta^k_m \\
&ii)& \;\;\; \sum\limits_m|f_m(z)|\leq M <\infty \;\; for \;any \; z
\in U
\end{eqnarray*}
(the symbol $\delta^k_m$ is equal to 1 if $m=k$ and 0 otherwise)

\end{Lemma}

\begin{proof}
By proposition 2 and lemma 12 in \cite{Toe} .
\end{proof}

\begin{Lemma} Let $\Omega$ be a domain in $C^n$, $f\in H(\Omega)$. If a
compact set $K$ and its neighborhood $G$ satisfy $K\subset
 G\subset\subset\Omega$ and $\rho=dist(K,\partial G)>0$, then
$$\sup\limits_{z\in K}|\frac{\partial
f}{\partial z_{j}}(z)|\leq\frac{\sqrt{n}}{\rho}\sup\limits_{z\in
G}|f(z)|.$$\end{Lemma}

\begin{proof}\hspace{4mm} Since $\rho=dist(K,\partial G)>0,$ for any
$a\in K$, the polydisc
$$
P_a=\left\{(z_1, \cdots, z_n)\in C^n: |z_j-a_j|
<\displaystyle\frac{\rho}{\sqrt{n}}, j=1,\cdots,n\right\}
$$
is contained in $G$. Using Cauchy inequality, we have
$$
\left|\displaystyle\frac{\partial f}{\partial z_j}(a)
\right|\leq\displaystyle\frac{\sqrt{n}}{\rho}
\sup\limits_{z\in\partial_0 P_a}|f(z)|\leq
\displaystyle\frac{\sqrt{n}}{\rho}\sup\limits_{z\in G}|f(z)|.
$$
So the Lemma follows.\end{proof}

\begin{Lemma} For fixed $0<\delta<1$, let $F_{\delta}=\{z\in U^n : \max\limits_{1\leq j\leq n}|z_j|>1-\delta \}.$ Then
$$\lim\limits_{r\rightarrow 1}\sup\limits_{||f||_\infty=1}\sup\limits_{z\in F_{\delta}^c}|f(z)-f(rz)|=0$$
for any $f\in H^{\infty}(U^n)$.\end{Lemma}

\begin{proof}
\begin{eqnarray*}
&&\sup\limits_{z\in F_{\delta}^c}|f(z)-f(rz)|\\
&=&\sup\limits_{z\in
F_{\delta}^c}|\sum\limits^{n}_{j=1}(f(rz_1,rz_2,\cdots,rz_{j-1},z_{j},\cdots,z_n)\\
&-&f(rz_1,rz_2,\cdots,rz_{j},z_{j+1},\cdots,z_n))|\\
&\leq& \sup\limits_{z\in
F_{\delta}^c}\sum\limits^{n}_{j=1}\left|\int^{1}_{r}|z_{j}\frac{\partial
f}{\partial
z_{j}}(rz_1,\cdots, rz_{j-1},tz_j,z_{j+1},\cdots,z^n)dt\right|\\
&\leq& (1-r)n\sup\limits_{z\in F_{\delta}^c}\left|\frac{\partial
f}{\partial z_{j}}(z)\right|.
\end{eqnarray*}
Consider $F_{\delta/2}^c$, then $F_{\delta}^c\subset F_{\delta/2}^c$
and $dist(F_{\delta/2}^c,\partial U^n)=\frac{\delta}{2}$.

From Lemma 3, we have
$$\sup\limits_{z\in F_{\delta}^c}\left|\frac{\partial f}{\partial z_{j}}(z)\right| \leq
\frac{2\sqrt{n}}{\delta}\sup\limits_{z\in F_{\delta/2}^c}|f(z)|.$$
 then
$$\sup\limits_{z\in
F_{\delta}^c}|f(z)-f(rz)|\leq\frac{2(1-r)n\sqrt{n}}{\delta}||f||_{\infty}.$$

Let $r\rightarrow 1$, the conclusion follows.\end{proof}

\section{Main theorem}

{\bf Theorem.}\hspace*{4mm} Let $\varphi,\psi:U^n \rightarrow U^n$ ,
and $C_{\varphi}-C_{\psi}:H^{\infty}(U^n)\rightarrow
H^{\infty}(U^n)$. Then
\begin{eqnarray*}&&\lim\limits_{\delta \rightarrow
0}\sup\limits_{z\in
  E_\delta}\max\limits_{1\leq j\leq
n}\beta(\varphi_j(z),\psi_j(z))\leq
||C_{\varphi}-C_{\psi}||_e\\
&\leq&\frac{4-4\sqrt{1-\lim\limits_{\delta \rightarrow
0}\sup\limits_{z\in E_{\delta}}\max\limits_{1\leq j \leq
n}\beta(\varphi_j(z),\psi_j(z))^2}}{\lim\limits_{\delta \rightarrow
0}\sup\limits_{z\in E_{\delta}}\max\limits_{1\leq j \leq
n}\beta(\varphi_j(z),\psi_j(z))}.\end{eqnarray*}

\begin{proof} We consider the upper estimate first.
 For fixed $0<r<1$, it is easy to check that both $C_{r\varphi}$
 and $C_{r\psi}$ are compact operators. Therefore
$$
||C_{\varphi}-C_{\psi}||_e\leq
||C_{\varphi}-C_{\psi}-C_{r\varphi}+C_{r\psi}||.$$
Now for any
$0<\delta<1$
\begin{eqnarray*}
&&\|C_{\varphi}-C_{\psi}-C_{r\varphi}+C_{r\psi}\|\\
&=&\sup\limits_{||f||_{\infty=1}}||(C_{\varphi}-C_{\psi}-C_{r\varphi}+C_{r\psi})f||_{\infty}\\
&=&\sup\limits_{||f||_{\infty}=1}\sup\limits_{z\in
U^n}|f(\varphi(z))-f(\psi(z))-f(r\varphi(z))+f(r\psi(z))|\\
&\leq& \sup\limits_{||f||_{\infty}=1}\sup\limits_{z\in
E_\delta}|f(\varphi(z))-f(\psi(z))-f(r\varphi(z))+f(r\psi(z))|\\
&+&\sup\limits_{||f||_{\infty}=1}\sup\limits_{z\in
E_\delta^c}|f(\varphi(z))-f(\psi(z))-f(r\varphi(z))+f(r\psi(z))|.
\end{eqnarray*}

From Lemma 3, we can choose $r$ sufficiently close to $1$ such that
the second term of the right hand side is less than any given
$\varepsilon$, and denote the first term by $I$.

Using Schwartz-Pick lemma and the monotony of function
$f(x)=\frac{2-2\sqrt{1-x^2}}{x}$, Then
\begin{eqnarray*}
I &\leq& \sup\limits_{||f||_{\infty}=1} \sup\limits_{z \in
E_{\delta}}(|f(\varphi(z))-f(\psi(z))|+|-f(r\varphi(z))+f(r\psi(z))|)\\
&=&\sup\limits_{z\in
E_{\delta}}\sup\limits_{||f||_{\infty}=1}(|f(\varphi(z))-f(\psi(z))|+|-f(r\varphi(z))+f(r\psi(z))|)\\
&=&\sup\limits_{z\in
E_{\delta}}(d_{U^n}(\varphi(z),\psi(z))+d_{U^n}(r\varphi(z),r\psi(z)))\\
&\leq& 2\sup\limits_{z\in
E_{\delta}}\frac{2-2\sqrt{1-\max\limits_{1\leq j \leq n}\beta(\varphi_j(z),\psi_j(z))^2}}{\max\limits_{1\leq j \leq n}\beta(\varphi_j(z),\psi_j(z))}\\
&=& \frac{4-4\sqrt{1-\sup\limits_{z\in E_{\delta}}\max\limits_{1\leq
j \leq n}\beta(\varphi_j(z),\psi_j(z))^2}}{\sup\limits_{z\in
E_{\delta}}\max\limits_{1\leq j \leq
n}\beta(\varphi_j(z),\psi_j(z))}.
\end{eqnarray*}
Let $\delta \rightarrow 0$, the upper estimate follows.

Now we turn to the lower estimate.

 Define $a_j=\lim\limits_{\delta\rightarrow 0}\sup\limits_{z\in
 E^j_\delta}\beta(\varphi_j(z),\psi_j(z))$. If we set $\delta_m=\frac{1}{m}$, then $\delta_m \rightarrow 0$ as
$m\rightarrow \infty$, and there exists $ z_m \in
 E^j_{\delta_m}$ such that $\lim\limits_{m\rightarrow
 \infty}\beta(\varphi_j(z_m),\psi_j(z_m))=a_j$.

Without loss of generality, we can assume
$|\varphi_j(z_m)|\rightarrow
 1$. Let $w_m=\varphi_j(z_m)$, by lemma 1, we further assume the
 subsequence by $w_m$, there exists a number $M_j$ and a sequence of
 functions $f_m\in H^\infty(D)$ such that
\begin{eqnarray*}
&i)& \;\;\; f_m(w_k)=\delta^k_m \\
&ii)& \;\;\; \sum\limits_m|f_m(w)|\leq M_j <\infty \;\; for \;any \;
w \in U
\end{eqnarray*}

Now for any $z \in U^n$, we define
 $\tilde{f_m}(z):=f_m(z^j)$, where
 $z^j$ is the $j^{th}$ component of $z$, then
 $\sum\limits_m|\tilde{f_m}(z)|\leq M_j
 <\infty$.

 Next we claim that $\tilde{f_m}$ converge weakly to $0$. In
 fact, let $\lambda \in H^\infty(U^n)^\ast$. For any integer $N$, there
 exist some unimodular sequence $\alpha_m$ such that
\begin{eqnarray*}
\sum^N_{m=0}|\lambda\tilde{f_m}| &=& \sum^N_{m=0} \lambda\tilde{f_m}
\alpha_m=\lambda(\sum^N_{m=0}
\lambda\tilde{f_m})\\
&\leq& ||\lambda|| ||\sum^N_{m=0} \alpha_m \tilde{f_m}||_\infty \leq
||\lambda||M_j
\end{eqnarray*}
Thus $\lambda \tilde{f_m} \rightarrow 0$,that is $\tilde{f_m}$
converge weakly to $0$.

 Putting functions
 $$g_m(z)=\frac{\tilde{f_m}(z)}{M_j}\frac{z^j-\psi_j(z_m)}{1-\overline{\psi_j(z_m)}z^j}$$then
 $||g_m||_\infty \leq 1$ and $g_m(z)$ converge weakly to $0$. So for any compact operator $K$, we have
 $||Kg_m||_\infty \rightarrow
 0$.

Now  we have
\begin{eqnarray*}
||C_{\varphi}-C_{\psi}-K|| &\geq& \limsup\limits_{m\rightarrow \infty}||(C_{\varphi}-C_{\psi}-K)g_m||_{\infty}\\
&\geq&\limsup\limits_{m\rightarrow
\infty}(||(C_{\varphi}-C_{\psi})g_m||_{\infty}-||Kg_m||_\infty)\\
  &=&\limsup\limits_{m\rightarrow \infty}\sup\limits_{z\in
  U^n}|g_m(\varphi(z))-g_m(\psi(z))|\\
  &\geq& \frac{1}{M_j}\limsup\limits_{m\rightarrow \infty}\sup\limits_{z\in
  U^n}|\frac{\varphi_j(z)-\psi(z_m)}{1-\overline{\psi_j(z_m)}\varphi_j(z)}\tilde{f_m}(\varphi(z))\\
  & & -\frac{\psi_j(z)-\psi_j(z_m)}{1-\overline{\psi_j(z_m)}\psi(z)}\tilde{f_m}(\psi(z))|\\
  &\geq&\frac{1}{M_j}\limsup\limits_{m\rightarrow \infty}\frac{|\varphi(z_m)-\psi(z_m)|}{1-\overline{\psi(z_m)\varphi(z_m)}}=\frac{1}{M_j}a_j.
\end{eqnarray*}

For the case$|\psi_j(z)|\rightarrow 1$, a similar argument can get
the same conclusion except $M_j$ is substituted by  a new constant
$M'_j$. If we set $\tilde{M_j}=\max\{M_j,M'_j\}$. Then for $1\leq
j\leq n$, and so we get the following estimate
\begin{eqnarray*}
||C_\varphi-C_\psi||_e&\geq& \frac{1}{\tilde{M_j}}\max\limits_{1\leq
j\leq n}\lim\limits_{\delta \rightarrow 0}\sup\limits_{z\in
  E^j_\delta}\beta(\varphi_j(z),\psi_j(z))\\
  &\geq \frac{1}{M}&\lim\limits_{\delta \rightarrow 0}\sup\limits_{z\in
  E_\delta}\max\limits_{1\leq j\leq
n}\beta(\varphi_j(z),\psi_j(z)).
\end{eqnarray*}
where $M=\max\limits_{1\leq j\leq n}\tilde{M_j}$.
\end{proof}

{\bf Corollary.}\hspace*{4mm}$C_{\varphi}-C_{\psi}$ is compact if
and only if$$\lim\limits_{\delta \rightarrow 0}\sup\limits_{z\in
  E_\delta}\max\limits_{1\leq j\leq
n}\beta(\varphi_j(z),\psi_j(z))=0.$$

\begin{proof}
By the inequality $\frac{1-\sqrt{1-x^2}}{x} \leq x$ for any $0<x\leq
1$ and $T$ is compact if and only if $||T||_e=0$
\end{proof}

\textbf {Remark.}\hspace*{4mm}If for any $j$, we have
$||\varphi_j||_{\infty}<1$ and $||\psi_j||_{\infty}<1$, then
$E_\delta=\emptyset$ when $\delta $ is small enough, without loss of
generality, we set
$$\lim\limits_{\delta \rightarrow 0}\sup\limits_{z\in
  E_\delta}\max\limits_{1\leq j\leq
n}\beta(\varphi_j(z),\psi_j(z))=0.$$

\end{document}